\documentclass{amsart}
\usepackage{amssymb}
\usepackage{amscd}
\newtheorem{Theorem}{Theorem}[section]
\newtheorem{Lemma}[Theorem]{Lemma}
\newtheorem{Corollary}[Theorem]{Corollary}
\newtheorem{Proposition}[Theorem]{Proposition}
\theoremstyle{definition}
\newtheorem{Definition}[Theorem]{Definition}
\newtheorem{Example}[Theorem]{Example}
\theoremstyle{remark}
\newtheorem{Remark}{Remark}

\font\sy=cmsy10
\font\ym=msbm10  
\newcommand{\Ker}{\hbox{\rm Ker}}
\newcommand{\End}{\rm End}
\newcommand{\Hom}{\rm Hom}
\newcommand{\Mat}{\rm Mat}
\newcommand{\C}{{\text{\ym C}}}
\newcommand{\N}{{\text{\ym N}}}
\newcommand{\cA}{{\hbox{\sy A}}}
\newcommand{\cB}{{\hbox{\sy B}}}
\newcommand{\cC}{{\hbox{\sy C}}}
\newcommand{\cD}{{\hbox{\sy D}}}
\newcommand{\cE}{{\hbox{\sy E}}}
\newcommand{\cF}{{\hbox{\sy F}}}
\newcommand{\cH}{{\hbox{\sy H}}}
\newcommand{\cM}{{\hbox{\sy M}}}

\newcommand{\cP}{{\hbox{\sy P}}}
\newcommand{\cR}{{\hbox{\sy R}}}
\newcommand{\cS}{{\hbox{\sy S}}}
\newcommand{\cT}{{\hbox{\sy T}}}
\newcommand{\cU}{{\hbox{\sy U}}}
\newcommand{\cV}{{\hbox{\sy V}}}

\title[Operator Categories]
{Notes on Operator Categories}
\author[Yamagami Shigeru]{Shigeru Yamagami}
\begin{document}
\maketitle
\begin{center}
Department of Mathematical Sciences\\
Ibaraki University\\
Mito, 310-8512, JAPAN\\
http://suuri.sci.ibaraki.ac.jp/\~{}yamagami/
\end{center} 

\begin{abstract}
The bicategory of normal functors between W*-categories is 
monoidally equivalent to the bicategory of W*-bimodules. 
\end{abstract}

\section*{Introduction}

Related to subfactor theory two kinds of tensor categories are used 
in describing quantum symmetry (\cite{L,O}): 
Given a von Neumann algebra $M$, 
we have the tensor category ${}_M\cB imod_M$ of 
$M$-$M$ bimodules on the one hand and 
the tensor category $\mathcal{E}nd(M)$ of endomorphisms on the other hand.
In spite of different appearances, they have a close similarity: 
For an endomorphism $\rho$ of $M$, assign the bimodule $L^2(M)\rho$, 
where $L^2(M)\rho$ is the standard Hilbert space of $M$ with the right 
action modified by $\rho$. Then we have canonical isomorphisms 
$L^2(M)\rho\otimes_M L^2(M)\sigma \to L^2(M)(\rho\circ \sigma)$ 
for endomorphisms $\rho$ and $\sigma$, which turn out to 
satisfy the condition of multiplicativity for monoidal functors 
and defines a fully faithful imbedding of the opposite of $\cE nd(M)$ 
into the tensor category ${}_M\cB imod_M$. 

If $M$ is of an infinite type, this monoidal imbedding gives an equivalence 
of categories, i.e., they contain the same information as categorical data. 
For a finite von Neumann algebra, however, the imbedding is not surjective 
because of the independant sizes of left and right modules in that case. 

Furthermore, the catgory $\cB imod$ of bimodules is more flexible than 
$\cE nd(M)$ in the sense that we can choose different algebras 
for left and right actions. The bimodules, 
together with the associative tensor products, then constitute 
a so-called bicategory 
(see \S 2 for explanations of the notion). 

With these observations in mind, 
the author has suspected the possibility to extend the tensor category 
$\cE nd(M)$ in such a way that the above imbedding can be extended to 
an equivalence of bicategories. 
The major purpose of this article is to give an affirmative answer 
to the above question in the framework of W*-categories: 
the bicategory of normal functors between W*-categories is 
monoidally equivalent to the bicategory $\cB imod$ of W*-bimodules. 

On the way of describing this result, we shall also review 
some of basic facts on operator categories 
(\cite{Ri}, \cite{GLR}, \cite{LR}), where 
we have fully used tensor products of W*-modules 
(the relative tensor product) to obtain concise expressions. 

Technically 
we have to use the modular theory in operator algebras 
on occasions, which is, however, algebraic (and formal) in nature rather than 
analytic if properly formulated (see \cite{Y1}).

As backgrounds of the subject, we refer, for example, to \cite{T} 
on operator algebras and to \cite{Mac, Pa} for generalities 
on category theory. 

\section{W*-categories}

By a linear category, we shall mean an essentially small 
category for which hom-sets 
are vector spaces over the complex number field $\C$ and 
the operation of taking compositions is complex-linear in variables involved. 

The restriction of essential smallness reflects our standing position 
that we shall not study operator algebras with the help of categorical 
languages but work with categories themselves of operator algebraic natures. 

A functor between linear categories is said to be linear if 
the operation on morphisms is linear. 
Recall that a functor is said to be faithful if it is injective on 
hom-sets (no commitment to objects). 

A linear category is called a {\bf *-category} if it is furnished with 
conjugate-linear involutions on hom-sets satisfying (i) 
$f^*: Y \to X$ for $f: X \to Y$ and (ii) $(g\circ f)^* = f^*\circ g^*$ 
for $f: X \to Y$ and $g: Y \to Z$. 

A *-category is called a {\bf C*-category} if hom-sets are Banach spaces 
such that (i) $\|g\circ f\| \leq \| g\|\,\| f\|$
for morphisms of the form $f:X \to Y$ and $g: Y \to Z$, 
(ii) $\| f^*f\| = \| f\|^2$ and (iii) $f^*f \geq 0$ for any morphism 
$f$. 

Note here that, by the conditions (i) and (ii), each $\End(X)$ is 
a C*-algebra and the meaning of positivity in (iii) is that for 
C*-algebras. 

A C*-category $\cC$ is called a {\bf W*-category} if each Banach space 
$\Hom(X,Y)$ is the dual of a Banach space. 
`Preduals' are uniquely determined by the C*-category $\cC$ 
as soon be checked (an analogue of Sakai's characterization, \cite{S}). 

A typical example of W*-categories is the category $\cR ep(A)$ 
of *-representations of a C*-algebra $A$ in
Hilbert spaces of a specified class 
with hom-sets given by intertwiners of representations. 

When $A$ is separable (i.e., having a countable norm-dense subset), 
we can define the much smaller W*-category
$\cS\cR ep(A)$ 
of *-representations of $A$ in separable 
Hilbert spaces of a specified class 
with hom-sets given by intertwiners of representations. 
If $A = \C$, $\cR ep(A)$ (resp.~$\cS\cR ep(A)$) 
is the category of Hilbert spaces (resp.~separable Hilbert spaces) 
of a specified class 
$\cH ilb$ (resp.~$\cS\cH ilb$) 
whose morphisms are bounded linear operators. 

Another typical example of W*-category is the category $\cP(M)$ 
of projections in a W*-algebra $M$: objects of $\cP(M)$ are projections in 
$M$ with hom-sets given by $\Hom(e,f) = fMe$ for projections $e$ and 
$f$ in $M$. 

A typical example of C*-category is the category $A$-$\cM od$ of 
Hilbert $A$-modules (again in a specified class) 
with $A$ a C*-algebra. 

Given *-categories $\cC$ and $\cD$, a functor $F: \cC \to \cD$ 
is called a {\bf *-functor} if $F$ is linear in morphisms and preserves 
the *-operation. 
If both of $\cC$ and $\cD$ are C*-categories, then 
a *-functor $F: \cC \to \cD$ is norm-decreasing: 
\[
\| F(f)\|^2 = \| F(f^*f)\| \leq \| f^*f\| = \| f\|^2
\]
for a morphism $f: X \to Y$ in $\cC$. 
The kernel of $F$ (which is a C*-subcategory of $\cC$) 
is then an analogue of closed *-ideals 
in C*-algebras and we have the exact sequence 
\[
0 \to \Ker(X,Y) \to \Hom(X,Y) \to F(\Hom(X,Y)) \to 0
\]
of C*-categories. 

Given *-categories $\cC$ and $\cD$, we define the category 
$\cH om(\cC,\cD)$ with objects given by *-functors and morphisms 
consisting of natural transformations. Recall here that we have 
stuck to essentially small categories, which enables us to keep 
$\cH om(\cC,\cD)$ essentially small. 
When $\cD$ is a C*-category, morphisms are restricted to 
be bounded so that 
the category $\cH om(\cC,\cD)$ of functors is again a C*-category. 

Given an object $X$ in a C*-category $\cC$ and a positive linear 
functional $\varphi$ on the C*-algebra $\End(X)$, we introduce 
the *-functor $F_\varphi: \cC \to \cH ilb$ by the GNS-construction: 
For an object $Y$ in $\cC$, the Hilbert space $F_\varphi(Y)$ is 
the completion of the vector space $\Hom(X,Y)$ with respect to 
the positive-semidefinite inner product 
\[
(y|y') = \varphi(y^*y')
\quad
\text{for $y,y': X \to Y$.}
\]
We shall often use the notation $y\varphi^{1/2}$ to distinguish 
the morphism $y$ in $\cC$ with the associated element in the Hilbert 
space $F_\varphi(Y)$. 

For a morphism $f: Y \to Z$ in $\cC$, $F_\varphi(f): F_\varphi(Y) \to F_\varphi(Z)$ 
is the bounded linear map defined by 
\[
F_\varphi(f)(y\varphi^{1/2}) = (fy)\varphi^{1/2}.
\]
These are well-defined by the positivity and the C*-norm assumption 
in the definition of C*-categories. 

Note that, if $X \cong Y$ and positive linear functionals 
$\varphi: \End(X) \to \C$, $\psi: \End(Y) \to \C$ are related by 
a *-isomorphism between $\End(X)$ and $\End(Y)$, then 
the functors $F_\varphi$ and $F_\psi$ are unitarily equivalent. 

Since (infinite) direct sums are permitted in the C*-category $\cH ilb$, 
we can define the *-functor $F: \cC \to \cH ilb$ as the direct sum 
of $\{ F_\varphi\}$, where $\varphi \in \End(X)^*_+$ with the family 
$\{ X\}$ representing the isomorphism classes of objects in $\cC$. 
Then $F$ is faithful on each C*-algebra $\End(X)$ and we have 
\[
\| F(f)\|^2 = \| F(f^*f)\| = \| f^*f\| = \| f\|
\]
for a morphism $f: X \to Y$, which implies that $F$ is isometric 
on hom-sets. Thus, given a finite family $\{ X_i\}_{1 \leq i \leq n}$ 
of objects in the C*-category $\cC$, the algebra of matrix form 
\[
\begin{pmatrix}
F(\Hom(X_1,X_1)) & \dots & F(\Hom(X_n,X_1))\\
\vdots &\ddots & \vdots\\
F(\Hom(X_1,X_n)) & \dots & F(\Hom(X_n,X_n))
\end{pmatrix}
\]
is a C*-algebra on the Hilbert space 
\[
\bigoplus_{i=1}^n F(X_i), 
\]
which is *-isomorphic to the *-algebra
\[
\begin{pmatrix}
\Hom(X_1,X_1) & \dots & \Hom(X_n,X_1)\\
\vdots &\ddots & \vdots\\
\Hom(X_1,X_n) & \dots & \Hom(X_n,X_n)
\end{pmatrix}. 
\]

\begin{Proposition}
Given a C*-category $\cC$ and a finite family $\{ X_i\}_{1 \leq i \leq n}$ 
of objects in $\cC$, the vector space 
\[
\begin{pmatrix}
\Hom(X_1,X_1) & \dots & \Hom(X_n,X_1)\\
\vdots &\ddots & \vdots\\
\Hom(X_1,X_n) & \dots & \Hom(X_nX_n)
\end{pmatrix}
\]
is a C*-algebra with each matrix component isometrically identified 
with the Banach spaces $\Hom(X_i,X_j)$. 
\end{Proposition} 

The above proposition is used to enlarge a C*-category $\cC$ 
so that it allows finite direct sums; 
consider the category $\widehat\cC$ for which 
objects are finite families of objects in $\cC$ and hom-sets are 
given by 
\[
\Hom(\{ X_i\}, \{ Y_j\}) 
= 
\begin{pmatrix}
\Hom(X_1,Y_1) & \dots & \Hom(X_m,Y_1)\\
\vdots & \ddots & \vdots\\
\Hom(X_1,Y_n) & \dots & \Hom(X_m,Y_n)
\end{pmatrix}
\]
with the norm induced from a bigger C*-algebra of matrix form. 

\begin{Corollary}
A C*-category is a W*-category if and only if 
\[
\begin{pmatrix}
\End(X) & \Hom(Y,X)\\
\Hom(X,Y) & \End(Y)
\end{pmatrix}
\]
is a W*-algebra for any pair $(X,Y)$ of objects in $\cC$. 
\end{Corollary}

Let $\cC$ be a W*-category. Since $\Hom(X,Y)$ is realized 
as a corner of a W*-algebra, we can define their $L^p$-extensions 
(see \cite{Ha} for $L^p$-theory on von Neumann algebras, cf.~also 
\cite{Y1,Y2}): 
the Banach space $L^p(X,Y)$ is defined to be the upper-right corner of 
the $L^p$-space associated to the W*-algebra in the above corollary. 
We then have the bounded bilinear map $L^q(X,Y)\times L^p(Y,Z) \to L^r(X,Z)$ 
if $ 1/r = 1/p + 1/q$ and $L^\infty(X,Y)$ is identified with 
$\Hom(Y,X)$. We also have the duality $L^p(X,Y)^* = L^q(Y,X)$ 
with $1/p + 1/q = 1$ and the predual of $\Hom(X,Y)$ is identified with 
the Banach space $L^1(X,Y)$. 

\begin{Proposition}
Given a finite family $\{ X_i\}$ of objects in a W*-category, 
the matrix algebra 
\[
\begin{pmatrix}
\Hom(X_1,X_1) & \dots & \Hom(X_n,X_1)\\
\vdots & \ddots & \vdots\\
\Hom(X_1,X_n) & \dots & \Hom(X_n,X_n)
\end{pmatrix}
\]
is a W*-algebra. 
Moreover, any W*-category is (in a unique way) 
extended so that it admits finite direct sums. 
\end{Proposition}

Let $\cC$ be a *-category and $\cD$ be a W*-category. 
Then $\cH om(\cC,\cD)$ is again a W*-category ($\cC$ being assumed to 
be essentially small and natural transformations being restricted 
to be bounded). 
For the identity functor $\text{id}_\cC$, $\End(\text{id}_\cC)$ is 
then a commutative W*-algebra $Z(\cC)$ given by 
\[
\{ t = \{ t_X\}; t_X \in \End(X), 
f t_X = t_Yf\ \text{for any $f \in \Hom(X,Y)$ in $\cC$}\}.
\]
We call $Z(\cC)$ the {\bf center} of $\cC$. 
If $\cC$ is a full subcategory of $M$-$\cM od$ including a faithful 
representation, then 
$Z(\cC)$ is naturally isomorphic to the center of $M$. 

For a morphism $f:X \to Y$ in $\cC$, the {\bf central support} of $f$ 
is defined to be the projection $c(f)$ in $Z(\cC)$ which is minimal 
among the projections $c$ in $Z(\cC)$ satisfying 
$c_Y f = f = f c_X$. 
For a family $\{ f_i\}$ of morphisms, we define its central support 
as the smallest projection in $Z(\cC)$ majorizing all 
$c(f_i)$'s: $c = \bigvee_i c(f_i)$. 
For a family $\{ X_i\}$ of objects, its central support is defined to 
be the central support of the family $\{ 1_{X_i} \}$. 

The central support $c$ of $f$ is calculated by 
\[
c_Z = \bigvee_{g: Z \to X} s(fg), 
\]
where $s(fg)$ denotes the support projection of the element 
$g^*f^*fg$ in $\End(Z)$. In particular, if $\Hom(X,Y) = \{ 0\}$, i.e., 
$X$ and $Y$ are disjoint, then $c(1_X) c(1_Y) = 0$. 

A family $\{ U_i\}$ of objects in a *-category $\cC$ is said to be 
{\bf generating} if 
the associated family of hom-functors 
$\Hom(U_i,\cdot): \cC \to \cV ec$ is faithful, i.e., 
the algebraic direct sum 
\[
\bigoplus_i \Hom(U_i,\cdot): \cC \to \cV ec
\]
is faithful. 
A single object $U$ in $\cC$ is called a {\bf generator} if 
the one-object family $\{ U\}$ is generating, i.e., 
the hom-functor $\Hom(U,\cdot): \cC \to \cV ec$ is faithful. 

\begin{Example}
Any representative family of isomorphism classes of objects is 
generating. 
\end{Example}

\begin{Lemma}
Let $\{ U_i\}$ be a family of objects in a W*-category $\cC$. 
Then the following conditions are equivalent. 
\begin{enumerate}
\item 
The family $\{ U_i\}$ is generating.
\item 
For any object $X$ and any projection $0 \not= p \in \End(X)$, 
we can find a $U_i$ such that 
$\Hom(U_i,pX) = p\Hom(U_i,X) \not= \{ 0\}$. 
\item 
The central support of $\{ U_i\}$ is equall to the family 
$\{ 1_X\}$ of identity morphisms. 
\item 
For any object $X$ in $\cC$, we can find a family of partial isometries 
$\{ u_{i,j}: U_i \to X\}$ such that 
\[
\sum_{i,j} u_{i,j} u_{i,j}^* = 1_X. 
\]
\end{enumerate}
\end{Lemma}

\begin{proof}
(i) $\Rightarrow$ (ii). If there exist an object $X$ and 
a projection $0 \not= p \in \End(X)$ such that 
$\Hom(U_i, pX) = \{ 0\}$ for any $i$, 
then $p$ can not be distinguished with $0$ under the hom-functor 
$\Hom(U_i,\cdot)$, which contradicts with the condition (i). 

(ii) $\Rightarrow$ (iii). If the central support $\{ c_X\}$ of 
$\{ U_i\}$ is different from $\{ 1_X\}$, we can find an object $X$ 
such that $c_X \not= 1_X$, which satisfies 
$(1_X - c_X)f = 0$ for any $i$ and any $f: U_i \to X$. 

(iii) $\Rightarrow$ (iv). 
The formula for the central support together with 
polar decompositions shows that for any object $X$, we can find 
$i$ and a non-zero partial isometry $u: U \to X$. 
Now the maximality argument gives the result. 

(iv) $\Rightarrow$ (i) is obvious.

\end{proof}

\begin{Lemma}[{\cite[Lemma~2.1]{Ro}}]
Let $X$, $Y$ and $Z$ be objects in $\cC$. 
Let $T: L^2(Y,X) \to L^2(Z,X)$ be a bounded linear map 
satisfying $T(\xi x) = T(\xi)x$ for $\xi \in L^2(Y,X)$ and 
$x \in \End(X)$. Then we can find an element $y \in \Hom(Y,Z)$ 
such that $T(\xi) = y\xi$ for $\xi \in L^2(Y,X)$. 
\end{Lemma}

\begin{proof}
Let $M$ be the W*-algebra associated to the family $\{ X, Y, Z\}$ 
in Proposition~1.3 with
$e$, $f$ and $g$ projections of $M$ to the component 
$X$, $Y$ and $Z$ respectively. 
Consider a bounded operator $T$ on $L^2(M)$ satisfying 
$T(\xi) = gT(f\xi e)e$ for $\xi \in L^2(M)$ and 
$T(\xi exe) = T(\xi)exe$ for $x \in M$. 

Since the commutant of a reduction is the induction of the commutant, 
the second condition on $T$ implies we can find $y \in M$ satisfying 
$T(\xi) = y\xi$. Now the restrictions on $T$ for the domain and range 
reveals that we can replace $y$ by the element $gyf$ and we are done.
\end{proof}

\begin{Proposition}
Given a generating family $\{ U_i\}$ in a W*-category $\cC$, 
let $M$ be the opposite of the von Neumann algebra 
$\bigoplus_{i,j} \Hom(U_i,U_j)$. Then the functor 
$F: \cC \to \cH ilb$ defined by 
\[
F(X) = \bigoplus_i L^2(X,U_i)
\]
with the obvious right action of $M^\circ$ on $F(X)$ 
gives a fully faithful imbedding of $\cC$ into the W*-category 
$M$-$\cM od$. 
\end{Proposition}

\begin{proof}
The surjectivity of the functor maps on morphisms is a consequence 
of the previous lemma because local solutions admit a convergent net 
by the compactness argument. 

On the other hand, if $f: X \to Y$ vanishes on the Hilbert space 
$F(X)$, then $f\Hom(U_i,X) \varphi_i^{1/2} = \{ 0\}$ for any $i$ and 
any $\varphi_i \in \End(U_i)^+_*$, whence 
$f\Hom(U_i,X) = \{ 0\}$ for $i$. 
Since $\{ U_i\}$ is generating, this implies $f = 0$. 
\end{proof}

\begin{Lemma}
For a von Neumann algebra $M$, the following are equivalent. 
\begin{enumerate}
\item 
$M$ has the separable predual.
\item 
$M$ has a faithful normal representation on a separable Hilbert 
space. 
\item 
The standard space $L^2(M)$ is separable. 
\end{enumerate}
\end{Lemma} 

\begin{proof}
Non-trivial is (i) $\Rightarrow$ (iii). 

Let $\{ f_i\}_{i \geq 1}$ be a countable dense set of $M_*$. 
Then we can find a family $\{ x_i \in M\}_{i \geq 1}$ such that 
$\| x_i \| \leq 1$ and $| f_i(x_i) | \geq \| f_i\|/2$, 
which turns out to be total in $M$ 
with respect to the weak*-topology. 
In fact, if not, we can find 
an element $x \in M$ and a functional $f \in M_*$ satisfying 
$f(x) = 1$ and $f(x_i) = 0$ for $i \geq 1$. For each integer $n \geq 1$, 
choose $i_n \in \N$ so that $\| f_{i_n} - f \| \leq 1/n$. 
Then the inequality 
\[
\frac{1}{2} \| f_{i_n} \| \leq 
|f_{i_n}(x_{i_n})| = | f_{i_n}(x_{i_n}) - f(x_{i_n}) | 
\leq \| f_{i_n} - f\| \leq \frac{1}{n}
\]
implies $\| f_{i_n} \| \to 0$ and hence 
$f = \lim_n f_{i_n} = 0$, which contradicts the choice $f(x) = 1$. 


The separability of the predual $M_*$ also ensures that 
we can find a faithful positive functional $\varphi$ in $M_*$. 
Now the set $\{ x_i \varphi^{1/2} \}_{i \geq 1}$ is total in $L^2(M)$ 
by the Kaplansky's density theorem and we are done. 
\end{proof}

\begin{Remark}~
\begin{enumerate}
\item
A similar argument shows that, if a Banach space $X$ has the separable 
dual Banach space, then $X$ itself is separable. 
Thus, the separabilty of $L^p(M)$ for some $1 < p < +\infty$ implies 
the separabilty of $L^q(M)$ with $1/p + 1/q = 1$ and then 
the predual is separable as a continuous image of 
$L^p(M)\times L^q(M)$. 
\item 
The argument in the above proof reveals that a von Neumann algebra 
is countably generated if it has the separable predual. 
The converse implication is, however, not true; 
the dual Banach space of the separable C*-algebra $C[0,1]$ is 
identified with the space $M[0,1]$ of complex measures in the interval 
$[0,1]$, which is not separable because $\| \delta_s - \delta_t\| = 2$ 
if $s \not= t$ in $[0,1]$. Then the double dual von Neumann algebra 
$M[0,1]^*$ is countably generated with the non-separable predual 
$M[0,1]$. 
\end{enumerate}
\end{Remark}

\begin{Definition}
A W*-category is said to be {\bf locally separable} 
if each $\Hom(X,Y)$ has the separable predual. 
A locally separable W*-category is said to be {\bf separable} 
if it admits a generating family indexed by a countable set. 
\end{Definition}

Given a W*-algebra $M$, 
we denote by $\cR ep(M)$ the W*-category of normal 
*-representations of $M$ on Hilbert spaces in a specified 
class. 
When $M$ has the separable predual, 
we denote by $\cS\cR ep(M)$ the W*-category of normal 
*-representations of $M$ on separable Hilbert spaces in a specified 
class. 

\begin{Proposition}~ 
\begin{enumerate}
\item 
A W*-category $\cC$ is equivalent to $\cR ep(M)$ for some W*-algebra $M$ 
if and only if $\cC$ admits subobjects 
(projections are associated to subobjects) and direct sums 
indexed by sets in a specified class. 
\item 
A W*-category $\cC$ is equivalent to $\cS\cR ep(M)$ for some W*-algebra $M$ 
of separable predual if and only if $\cC$ is separable and 
admits countably many direct 
sums as well as subobjects. 
\end{enumerate}
\end{Proposition}

\begin{proof}
Clearly the condition is necessary. 
On the other hand, we have the fully faithful imbedding 
$F:\cC \to \text{$M$-$\cM od$}$ by Proposition~1.7 and any $M$-module is 
isomorphic to $F(X)$ for some object $X$ in $\cC$ by the condition 
(cf.~the structure theorem of normal representations \cite{D,T}). 
\end{proof}

For functors between W*-categories, it is reasonable to restrict to 
normal ones: a *-functor $F: \cC \to \cD$ of W*-categories 
is said to be {\bf normal} if $F$ is weak*-continuous on each 
$\Hom(X,Y)$. 

\begin{Lemma}[{\cite[Prop.~4.7]{Ri}}]
Let $\{ U_i\}$ be a generating family of objects in $\cC$. 
A *-functor $F: \cC \to \cD$ between W*-categories is normal if and only 
if it is normal on each W*-algebra $\End(U_i)$.
\end{Lemma}

\begin{proof}
Assume that $F$ is normal on $\End(U_i)$'s. 

We will first show that $F$ is normal on each $\End(X)$: 
for any $\varphi \in \End(F(X))_*^+$, the functional 
$\End(X) \ni x \to \varphi(F(x))$ is weak*-continuous. 
To see this, it suffices to check the $\sigma$-strong continuity on 
the unit ball of $\End(X)$ by a well-known result 
(\cite[Theorem I.3.1]{D}). 
Let $x_\alpha \to x$ in the $\sigma$-strong topology 
with $x_\alpha$ and $x$ in the unit ball. 
Here choose a family of partial isometries 
$\{ u_{i,j}: U_i \to X \}$ satisfying $\sum u_{i,j} u_{i,j}^* = 1_X$. 

Then $u_{i,j} (x_\alpha - x)^*(x_\alpha - x) u_{i,j}^* \in \End(U_i)$ 
converges to $0$ in the weak*-topology and hence by the normality 
of $F$ on $\End(U_i)$, we have 
\[
F(u_{i,j})^* F(x_\alpha - x)^* F(x_\alpha - x) F(u_{i,j}) 
= F(u_{i,j}^* (x_\alpha - x)(x_\alpha -x )u_{i,j}) \to 0
\]
in the weak*-topology, i.e., 
$F(x_\alpha - x) F(u_{i,j})$ converges to $0$ in the $\sigma$-strong topology 
for each index $(i,j)$. Thus
\[
(\varphi^{1/2}| F(x_\alpha - x) F(u_{i,j}u_{i,j}^*) \varphi^{1/2}) \to 0 
\quad 
\text{as $\alpha \to \infty$}
\]
for each $(i,j)$. 

Since $\| F(x_\alpha - x)\| \leq 2$ and 
\[
\sum (\varphi^{1/2}|(F(u_{i,j}u_{i,j}^*)\varphi^{1/2}) = (\varphi^{1/2}|\varphi^{1/2})
\]
is finite, the usual argument shows that 
$\varphi(F(x_\alpha - x)) = (\varphi^{1/2}|F(x_\alpha - x)\varphi^{1/2})$ converges to $0$, 
proving the normality of $F$ on $\End(X)$.  

Now let $f_\alpha$ be a net in 
$\Hom(X,Y)$ which converges to $0$ in $\sigma$*-strong topology. 
Then $f_\alpha^*f_\alpha \to 0$ in weak*-topology and hence 
$F(f_\alpha)^* F(f_\alpha) = F(f_\alpha^*f_\alpha) \to 0$, 
i.e., $F(f_\alpha) \to 0$ in 
$\sigma$-strong topology. 
Thus the *-homomorphism 
\[
\begin{pmatrix}
\End(X) & \Hom(Y,X)\\
\Hom(X,Y) & \End(Y)
\end{pmatrix}
\to 
\begin{pmatrix}
\End(F(X)) & \Hom(F(Y),F(X))\\
\Hom(F(X),F(Y)) & \End(F(Y))
\end{pmatrix}
\]
between W*-algebras induced from $F$ is $\sigma$-strongly continuous. 
Then it is weak*-continuous by the well-known result on topologies 
of W*-algebras (\cite[Theorem I.3.1]{D}) and we are done. 
\end{proof}

\begin{Corollary}
An equivalence of W*-categories is automatically normal.    
\end{Corollary}

\begin{proof}
This follows from the automatical normality of *-isomorphisms 
between W*-algebras (the nomality being the condition on 
the order structure of hermitian elements). 
\end{proof}

Let $M$ and $N$ be W*-algebras and $H$ be an $N$-$M$ bimodule. 
Then we can define a normal *-functor $\cR ep(M) \to \cR ep(N)$ by 
\[
F({}_M X) = {}_N H\otimes_M X 
\]
(see \cite{Sa} further information on relative tensor products, 
cf.~also \cite{Y1,Y2}). 

Conversely we have the following reformulation of Rieffel's 
theorem on Morita equivalences. 

\begin{Proposition}
Any normal *-functor $F: \cR ep(M) \to \cR ep(N)$ 
(resp.~ $F: \cS\cR ep(M) \to \cS\cR ep(N)$ with 
$M$ and $N$ having separable preduals) is unitarily 
equivalent to the one associated to an $N$-$M$ bimodule $H$ 
(resp.~a separable $N$-$M$ bimodule). 

The bimodule $H$ is uniquely determined up to unitary isomorphisms
by the functor $F$: 
$H = {}_N F({}_ML^2(M))_M$ with the right action of 
$M = \End({}_ML^2(M))^\circ$ on $H$ given by the normal *-homomorphism 
$\End({}_ML^2(M)) \to\End({}_NF({}_ML^2(M)))$ induced by 
the functor. 
\end{Proposition}

\begin{proof}
Let ${}_MX$ be an $M$-module. By the structure theorem of normal 
*-homomorphism (\cite[Theorem I.4.3]{D}), 
we can find an index set $I$ and a projection in 
$\Mat_I(M)$ so that 
${}_MX \cong {}_M L^2(M)^{\oplus I}p$. 
We shall then construct a unitary intertwiner 
\[
{}_N F({}_ML^2(M))\otimes_M X \to {}_N F({}_MX)
\]
by the commutativity of 
\[
\begin{CD}
H\otimes_MX @>>> F(X)\\
@VVV @VVV\\
H\otimes_M L^2(M)^{\oplus I}p @>>> 
F({}_ML^2(M)^{\oplus I}p).
\end{CD}
\]
Here the bottom line is a unitary map given by the composition of 
\begin{gather*}
{}_N F({}_ML^2(M))\otimes_M L^2(M)^{\oplus I}p 
= {}_N F({}_ML^2(M))^{\oplus I}p\\
= {}_N F({}_ML^2(M)^{\oplus I})F(p)
= {}_NF({}_ML^2(M)^{\oplus I}p). 
\end{gather*}

By the multiplicativity of $F$ on morphisms, the unitary map 
$H\otimes_M X \to F(X)$ 
is independent of the choice of an isomorphism 
${}_M X \to {}_M L^2(M)^{\oplus I}p$ and behaves naturally for 
intertwiners. 
\end{proof}

A W*-category is said to be of {\bf type I} if the matrix 
W*-algebra 
\[
\begin{pmatrix}
\End(X) & \Hom(Y,X)\\
\Hom(X,Y) & \End(Y)
\end{pmatrix}
\]
is of type I (in the sense of Murray-von Neumann) 
for any pair $(X,Y)$ of objects. 

The following is an easy reformulation of the structure theorem 
on W*-algebras of type I. 

\begin{Proposition}[{cf.~\cite[Theorem~8.10]{Ri}}]
Assume that a W*-category $\cC$ is saturated for taking subobjects. 
Then $\cC$ is of type I if and only if 
we can find a generating family $\{ U_i\}$ satisfying 
(i) $\End(U_i)$ is commutative and (ii) 
$\Hom(U_i,U_j) = \{ 0\}$ for $i \not= j$. 
\end{Proposition}

\section{W*-Tensor Categories}

Recall that a {\bf tensor category} is a linear category 
$\cT$ together with a functor $\Phi: \cT\times \cT \to \cT$ 
and a natural isomorphism (the associativity constraint) 
$a: \Phi(\Phi\times \text{id}_\cT) \to \Phi(\text{id}_\cT\times \Phi)$ 
satisfying the so-called pentagonal condition. 
The tensor product notation is often used to denote the functor 
$\Phi$: $\Phi(X,Y) = X\otimes Y$ and 
$\Phi(f,g) = f\otimes g$ for objects $X$, $Y$ and morphisms 
$f$, $g$ in $\cT$. 

When $\cT$ is a *-category and $\Phi$ preserves the *-operation, i.e., 
$(f\otimes g)^* = f^*\otimes g^*$ for morphisms $f$, $g$ in $\cT$, 
$\cT$ is called a {\bf *-tensor category}. 
A *-tensor category is called a {\bf C*-tensor category} if 
it is based on a C*-category 
and the associativity constraint is unitary. 

A {\bf W*-tensor category} is, by definition, a C*-tensor category $\cT$ 
with the tensor product functor $\Phi: \cT\times \cT \to \cT$ is 
{\bf binormal} in the sense that it is 
separately normal on each variables. 

Similar adjective definitions work for bicategories; 
we can talk about *-bicategories, C*-bicategories and W*-bicategories. 

Recall that a bicategory consists of labels, $A,B,C, \dots$, 
categories $\cH om(A,B)$ indexed by pairs of labels, 
functors 
\[
\Phi_{A,B,C}: \cH om(B,A) \times \cH om(C,B) \to \cH om(C,A)
\]
indexed by triplets of labels together with natural isomorphisms 
\[
a_{A,B,C,D}: \Phi_{A,C,D} (\Phi_{A,B,C}\times \text{id}_{\cH om(D,C)}) 
\to \Phi_{A,B,D} (\text{id}_{\cH om(B,A)}\times \Phi_{B,C,D}) 
\]
indexed by quadruplets of labels and satisfying 
the pentagonal relation.

The functor $\Phi_{A,B,C}$ is often denoted by the notation of 
composition, which reflects the view-point that a bicategory 
is a `categorization' of hom-sets as hom-categories.

Here is another view-point from which we regard the hom-category 
$\cH om(A,B)$ as a resemblance of the category of $B$-$A$ bimodules 
(if labels represent algebras) with the notation 
$\cH om(A,B) = {}_B\cM_A$. 
Then the functor $\Phi_{A,B,C}$ is consequently denoted by 
the tensor product notation: For objects $X$ in ${}_A\cM_B$ 
and $Y$ in ${}_B\cM_C$, $\Phi_{A,B,C}(X,Y)$ is denoted by 
$X\otimes_B Y$. Similarly for morphisms. 

A typical example of W*-bicategory is provided by bimodules 
with normal actions of W*-algebras. 

A bicategory is said to be strict if the natural isomorphisms 
$a_{A,B,C,D}$ are identities, i.e., 
\[
\Phi_{A,C,D} (\Phi_{A,B,C}\times \text{id}_{\cH om(D,C)}) 
= \Phi_{A,B,D} (\text{id}_{\cH om(B,A)}\times \Phi_{B,C,D}) 
\]
and $a_{A,B,C,D}$ is the identity for each quadruplet. 

A typical example of strict bicategory is provided by categories of 
functors: 
Let $F,F': \cC \to \cD$ and $G, G': \cD \to \cE$ be functors. 
Then, given natural transformations $s: F \to F'$ and 
$t: G \to G'$, we can associate the natural transformation 
$G\circ F \to G'\circ F'$ by the commutative diagram
\[
\begin{CD}
G(F(X)) @>{G(s_X)}>> G(F'(X))\\
@V{t_{F(X)}}VV @VV{t_{F'(X)}}V\\
G'(F(X)) @>>{G'(s_X)}> G'(F'(X))\ .
\end{CD}
\]
Thus $F = {}_\cD F_\cC$, $G = {}_\cE G_\cD$ and 
${}_\cE G\otimes_\cD F_\cC = {}_\cE (G\circ F)_\cC$ 
with the unit objects given by identity functors. 
Note that 
the commutativity of the above diagram expresses the identity 
$(t\otimes_\cD 1_{F'}) (1_G \otimes_\cD s) 
= (1_{G'}\otimes_\cD s) (t\otimes_\cD 1_F)$. 
The (strict) associativity for tensor products of morphisms 
is also immediate. 

In particular, the category $\cE nd(\cC)$ of functors from $\cC$ into 
itself is a strict monoidal category. 

\begin{Example}
If $\cC$ is a one-object category, 
objects of $\cE nd(\cC)$ are endomorphisms of the algebra 
$A = \End(\text{pt})$ with hom-sets given by 
\[
\Hom(\rho,\sigma) = 
\{ a \in A; a \rho(x) = \sigma(x) a, \ 
\forall x \in A\}.
\]
The monoidal structure takes the form 
$\rho\otimes \sigma = \rho\circ \sigma$ for $\rho$, $\sigma \in \End(A)$ 
and $a\otimes b = a\rho(b) = \rho'(b)a$ 
for $a \in \Hom(\rho,\rho')$, $b \in \Hom(\sigma,\sigma')$. 

When $\cC$ is a W*-category, $A$ is a W*-algebra and 
objects in $\cE nd(\cC)$ are normal *-endomorphisms of $A$. 

The tensor category $\cE nd(\cC)$ is also denoted by $\cE nd(A)$ 
(see \cite{L} for more information on $\cE nd(A)$). 
\end{Example}

\begin{Lemma}[{\cite[Theorem 7.13]{GLR}, \cite[Lemma 2.1]{DPR}}]
Let $\cC$, $\cD$ be W*-categories and $\{ U_i\}$ be 
a generating family in $\cC$. 
Let $\cU$ be the full subcategory of $\cC$ consisting of objects 
in $\{ U_i\}$. 

Then the restriction (of functors and natural transformations) 
\[
\cH om(\cC,\cD) \ni F \mapsto F|_\cU \in \cH om(\cU,\cD)
\]
gives a fully faithful imbedding of W*-categories. 
Here $\cH om(\cC,\cD)$ and $\cH om(\cU,\cD)$ are 
W*-categories of normal *-functors and natural transformations. 

More concretely, given a natural transformation 
$\{ t_i: F(U_i) \to G(U_i)\}$ between normal *-functors 
$F|_\cU$ and $G|_\cU$, the natural transformation 
$t_X: F(X) \to G(X)$ is recovered by the formula 
\[
t_X = \sum_{i,j} F(u_{i,j}) t_i G(u_{i,j})^*, 
\]
where $\{ u_{i,j}: U_i \to X\}$ is a family of partial isometries 
satisfying $\sum_{i,j} u_{i,j} u_{i,j}^* = 1_X$. 
\end{Lemma}

\begin{proof}
Let $X$ be an object in $\cC$. 
Since $\{ U_i \}$ is a generating family, 
we can find a family $\{ u_{i,j}: U_i \to X\}$ 
of partial isometries such that $\sum_{i,j} u_{i,j}u_{i,j}^* = 1_X$. 

Given a natural transformation $t: F \to G$, we have 
$t_X F(u_{i,j}) = G(u_{i,j}) t_i$ for any $(i,j)$ and then 
\[
t_X = \sum_{i,j} t_X F(u_{ij}u_{i,j}^*) = 
\sum_{i,j} G(u_{i,j}) t_i F(u_{i,j})^*  
\]
by the normality of $F$. Thus the restriction is injective 
on natural transformations. 

Conversely, given a natural transformation 
$\{ t_i: F(U_i) \to G(U_i) \}$ between $F|_\cU$ and $G|_\cU$ and 
an object $X$ in $\cC$, write 
\[
t_X = \sum_{i,j} G(u_{i,j}) t_i F(u_{i,j})^*, 
\]
which is a morphism in $\Hom(F(X),G(X))$. 

If $Y$ is another object in $\cC$ with a family 
$\{ v_{k,l}: U_k \to Y \}$ of partial isometries 
satisfying $\sum_{k,l} v_{k,l} v_{k,l}^* = 1_Y$, we associate 
another morphism $t_Y: F(Y) \to G(Y)$. 
Then, for any morphism $f: X \to Y$ in $\cC$, we have 
\begin{align*}
G(f) t_X  &= \sum_{i,j} G(fu_{i,j}) t_i F(u_{i,j})^*\\
&= \sum_{i,j} \sum_{k,l} G(v_{k,l} v_{k,l}^* f u_{i,j}) 
t_i F(u_{i,j})^*\\
&= \sum_{i,j} \sum_{k,l} G(v_{k,l}) t_k 
F(v_{k,l}^* f u_{i,j}) F(u_{i,j})^*\\
&= t_Y F(f)
\end{align*}
by the naturality of $\{ t_i\}$ and the normality of $F$, $G$. 

If we take $f = 1_X$ ($X = Y$ particularly), then the above formula 
means that the morphism $t_X$ is well-defined. 
Thus, the natural rtansformation $\{ t_X \}$ is recovered from 
$\{ t_i\}$. 
\end{proof}

If we restrict ourselves to W*-categories, normal *-functors and 
bounded natural transformations, then we obtain 
the strict W*-bicategory $\cF unct$. 
Recall that the tensor category $\cE nd(M)$ in Example~2.1 is a part 
of $\cF unct$. 
In accordance with our separability notation, we denote by 
$\cS\cF unct$ the W*-bicategory of normal functors 
between separable W*-categories. 
Moreover we have the W*-bicategory $\cF unct_I$ of normal functors 
between W*-categories of type I. 

We shall now give a fully faithful imbedding of $\cF unct$ 
into the bicategory $\cB imod$ of W*-bimodules. 

To this end, we first enlarge relevant W*-categories so that they admit 
generators; we shall take a generator $U_\cA$ for each W*-category 
$\cA$ and define the W*-algebra $A$ by $A = \End(U_\cA)$ 
($B = \End(U_\cB)$ and son on). 
Then, by the imbedding theorem of W*-categories, 
we have fully faithful imbedding $\Phi_\cA$ of $\cA$ into the category of 
right $A$-modules: 
\[
\Phi_\cA: X \mapsto L^2(X,U_\cA)_A.
\]
To each functor $F: \cA \to \cB$, we associate the right $B$-module by 
\[
\Phi(F) = L^2(F(U_\cA),U_\cB). 
\]
Note here that the functor $F$ has the unique extension to 
the enlargements of $\cA$ and $\cB$ (cf.~Lemma~2.2). 
The Hilbert space $\Phi(F)$ admits the left action of $A$ by 
\[
a\xi = F(a)\xi, 
\quad 
a \in A,\ \xi \in L^2(F(U_\cA),U_\cB),
\]
which clearly commutes with the right action of $B$. 
Thus $\Phi(F)$ is an $A$-$B$ module. 

Moreover, given a natural transformation $t: F \to F'$ in 
$\cH om(\cA,\cB)$, the left multiplication of $t_{U_\cA} \in A$ 
defines an $A$-$B$ intertwiner between $\Phi(F)$ and $\Phi(F')$. 
In this way, we obtain a *-functor  
$\Phi: \cF unct \to \cB imod$. 
By the previous lemma, $\Phi$ is fully faithful. 

Given a functor $F: \cA \to \cB$ and an object $X$ in $\cA$, 
we define a linear map 
$m_{X,F}: L^2(X,U_\cA)\otimes_A \Phi(F) \to L^2(F(X),U_\cB)$ by 
\[
L^2(X,U_\cA)\otimes_A L^2(F(U_\cA),U_\cB) 
\ni x\varphi^{1/2}\otimes_{\varphi^{-1/2}} f\psi^{1/2} 
\mapsto F(x)f \psi^{1/2} 
\in L^2(F(U_\cA),U_\cB),
\]
which is clearly $B$-linear and isometric because of  
\[
(x\varphi^{1/2}\otimes_{\varphi^{-1/2}} f\psi^{1/2} | 
x\varphi^{1/2}\otimes_{\varphi^{-1/2}} f\psi^{1/2}) 
= (f\psi^{1/2} | F(x^*x) f\psi^{1/2}) 
= \| F(x)f \psi^{1/2} \|^2,
\]
where $\varphi \in A^+_*$ and $\psi \in B^+_*$. 

To see the surjectivity of $m_{X,F}$, 
choose a family $\{ x_i: U_\cA \to X \}$ 
of partial isometries satisfying $\sum_i x_ix_i^* = 1_X$. 
Then, for $y: U_\cB \to F(X)$, the normality of $F$ shows that 
\[
y = F\left(\sum_i x_i x_i^*\right)y = \sum_i F(x_i) F(x_i^*)y 
= \sum_i F(x_i) f_i, 
\]
where $f_i = F(x_i^*)y$ is in $\Hom(U_\cB, F(U_\cA))$. 
Thus the image of $m_{X,F}$ is dense in $L^2(F(X),U_\cB)$ 
and $m_{X,F}$ gives a unitary map. 

Let $G: \cB \to \cC$ be another normal *-functor. 
We can then define a unitary map 
$m_{F,G}: \Phi(F)\otimes_B \Phi(G) \to \Phi(G\circ F)$ by 
\[
m_{F(U_\cA),G}: 
L^2(F(U_\cA),U_\cB)\otimes_B L^2(G(U_\cB),U_\cC) 
\to L^2(GF(U_\cA),U_\cC), 
\]
which is $A$-$C$ linear by the definition of $m_{X,F}$ and the left action 
of $A$. The explicit form of these multiplication maps also shows 
the identity 
\[
m_{F(X),G}(m_{X,F}\otimes 1_{\Phi(G)}) = 
m_{X,GF}(1_{\Phi(X)}\otimes m_{F,G}). 
\]
In other words, the following diagram commutes. 
\[
\begin{CD}
\Phi(X)\otimes_A\Phi(F)\otimes_B\Phi(G) @>>> \Phi(X)\otimes_A \Phi(GF)\\
@VVV @VVV\\
\Phi(F(X))\otimes_B \Phi(G) @>>> \Phi(GF(X)). 
\end{CD}
\]

Summarizing the argument so far, we obtain the following. 

\begin{Theorem}
The opposite of the bicategory $\cF unct$ 
($\cS\cF unct$ or $\cF unct_I$ respectively) 
is monoidally equivalent to 
the bicategory $\cB imod$ of W*-bimodules 
($\cS\cB imod$ of separable W*-bimodules or $\cC\cB imod$ of 
W*-bimodules with actions of commutative W*-algebras respectively). 
\end{Theorem}

\begin{proof}
We have a fully faithful monoidal imbedding 
$\cF unct \to \cB imod$ by the above argument, whose image covers 
every bimodule up to unitary isomorphisms by Proposition~1.13. 
\end{proof}

\end{document}